\documentclass [11pt]{article}
\usepackage{graphicx}

\usepackage{amssymb}
\usepackage{amsmath}
\usepackage{enumerate}

\newcommand{\N}{\mathbb{N}}

\newcommand{\J}{\mathbb{J}}

\newtheorem{definition}{{\bf Definition}}[section]
\newtheorem{theorem}[definition]{{\bf Theorem}}

\newtheorem{corollary}[definition]{{\bf Corollary}}
\newtheorem{proposition}[definition]{\noindent {\bf Proposition}}

\newtheorem{example}[definition]{\noindent {\bf Example}}

\newtheorem{problem}[definition]{\noindent {\bf Problem}}
\newtheorem{fact}[definition]{\noindent {\bf Fact}}
\newtheorem{remark}[definition]{\noindent {\bf Remark}}

\newtheorem{lemma}[definition]{\noindent {\bf Lemma}}
\def\endproof{\hfill {\kern 6pt\penalty 500
\raise -0pt\hbox{\vrule \vbox to5pt {\hrule width 5pt
\vfill\hrule}\vrule}}}
\def\centerpicture #1 by #2 (#3){\leavevmode
        \vbox to #2{
        \hrule width #1 height 0pt depth 0pt
        \vfill
        \special{pictfile #3}}}
\begin{document}

\title{A characterization of well-founded algebraic lattices }
\author{Ilham Chakir\\
\small Math\'ematiques,\\
\small Universit\'e Hassan $1^{er}$,\\
\small Facult\'e des Sciences et Techniques,\\
\small Settat, Maroc\\
\small {\rm e-mail: ilham.chakir@univ-lyon1.fr} \and Maurice Pouzet \footnote{Supported by INTAS}\\
\small UFR de Math\'ematiques,\\
\small Universit\'e Claude-Bernard, \\
\small $43$, Bd. du $11$ Novembre $1918$,\\
\small $69622$ Villeurbanne, France\\
\small {\rm e-mail: pouzet@univ-lyon1.fr}}


\date{\today}
\maketitle
\begin{abstract}We characterize well-founded  algebraic lattices 
by means of forbidden subsemilattices of the join-semilattice made of their compact elements.  More specifically, we show that  an algebraic
 lattice $L$ 
 is well-founded if and only if $K(L)$, the join-semilattice of
compact elements of $L$, is well-founded and contains neither  $[\omega]^{<\omega}$, nor $\underline\Omega(\omega^*)$
 as a join-subsemilattice.   As an immediate corollary, we get that an  algebraic modular lattice $L$ is well-founded if and only if 
$K(L)$ is well-founded and contains no infinite independent set.  If $K(L)$ is a join-subsemilattice of $I_{<\omega}(Q)$, 
the set of finitely generated initial segments of a well-founded poset $Q$, then $L$ is well-founded if and only if $K(L)$ is well-quasi-ordered.\end{abstract}

\section{Introduction and synopsis of results}
Algebraic lattices \index{algebraic lattice} and join-semilattices
\index{join-semilattice} (with a 0) are two aspects of the same
thing, as expressed in the following basic result.
\begin{theorem} \cite{Hofmann}, \cite{grat} The collection $J(P)$ of ideals \index{ideal} of a join-semilattice $P$, once ordered by inclusion,
is an algebraic lattice
and the subposet $K(J(P))$ of its compact elements \index{compact
element} is isomorphic to $P$. Conversely, the subposet $K(L)$ of
compact elements of an algebraic lattice $L$ is a join-semilattice
with a $0$ and $J(K(L))$ is isomorphic to $L$.\end{theorem}

In this paper, we characterize well-founded algebraic lattices by
means of forbidden join-subsemilattices of the join-semilattice made
of their compact elements. In the sequel  $\omega$  denotes the
chain of non-negative integers, and when this causes no confusion,
the first infinite cardinal \index{cardinal} as well as  the first
infinite ordinal \index{ordinal}. We denote $\omega^*$ the chain of
negative integers.  We  recall that a poset $P$ is {\it
well-founded} \index{well-founded} provided that every non-empty
subset of $P$ has a minimal element. With the Axiom of dependent
choices, this amounts to the fact that $P$ contains no subset
isomorphic to $\omega^*$. Let $\Omega(\omega^*)$ be the set $[
\omega]^2$ of two-element subsets of $\omega$, identified to pairs
$(i,j)$, $i<j<\omega$, ordered so that $(i,j)\leq (i',j')$ if and
only if $i'\leq i$ and $j\leq j'$ w.r.t. the natural order on
$\omega$. Let $\underline\Omega(\omega^*):=\Omega(\omega^*) \cup
\{\emptyset \}$ be obtained by adding a least element. Note that
$\underline\Omega(\omega^*)$ is isomorphic to the set of bounded
intervals of $\omega$ (or $\omega^*$) ordered by inclusion. Moreover
$\underline\Omega(\omega^*)$ is a join-semilattice
($(i,j)\vee(i',j')=(i\wedge i',j\vee j')$). The join-semilattice
$\underline\Omega(\omega^*)$ embeds in $\Omega(\omega^*)$ as a
join-semilattice; the advantage of $\underline\Omega(\omega^*)$
w.r.t. our discussion is to have a zero. Let $\kappa$ be a cardinal
number, e.g. $\kappa:= \omega$; denote $[\kappa]^{<\omega}$  (resp.
$\mathfrak {P}(\kappa)$ ) the set, ordered by inclusion, consisting
of finite (resp. arbitrary ) subsets of $\kappa$. The posets
$\underline\Omega(\omega^*)$ and $[\kappa ]^{<\omega}$ are
well-founded lattices, whereas  the algebraic lattices $J(\underline
\Omega(\omega^*))$  and $J([\kappa]^{<\omega})$ ($\kappa$ infinite)
are not well-founded (and we may note that $J([\kappa]^{<\omega})$
is isomorphic to $\mathfrak {P}(\kappa)$). As  a poset
$\underline\Omega(\omega^*)$ is isomorphic to  a subset  of
$[\omega]^{<\omega}$, but not as a join-subsemilattice. This is our
first result.

\begin{proposition}\label{w-f}
$\underline\Omega(\omega^*)$ does not embed in
$[\omega]^{<\omega}$ as a join-subsemilattice;
 more generally, if $Q$ is a well-founded poset then $\underline \Omega(\omega^*)$ does
not embed as a join-subsemilattice into $I_{<\omega}(Q)$, the
join-semilattice made of finitely generated initial segments of
$Q$.
\end{proposition}
\begin{figure}
\begin{center}
\includegraphics[width=2.5in]{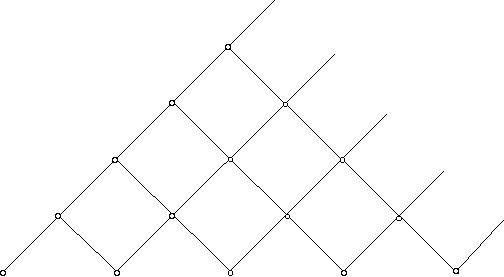}
\caption{$\Omega(\omega^*)$}
\end{center}
\label{fig:omega}
\end{figure}
Our next result expresses that $\underline
\Omega(\omega^*)$ and  $[\omega]^{<\omega}$ are unavoidable
examples of well-founded join-semilattices whose set of ideals is
not well-founded.
\begin{theorem}\label{thm4}
An algebraic lattice $L$  is well-founded if and only if $K(L)$ is
well-founded  and contains no join-subsemilattice  isomorphic to
$\underline \Omega(\omega^*)$ or to $[\omega]^{<\omega}$.\\
\end{theorem}

The fact that a join-semilattice $P$ contains a
join-subsemilattice isomorphic to $[\omega]^{<\omega}$  amounts to
the existence of an infinite independent set. Let us recall that a
subset $X$ of a join-semilattice $P$ is {\it independent}
\index{independent} if $x\not \leq \bigvee F$ for every $x\in X$
and every non-empty finite subset $F$ of $X\setminus \{x\}$.
Conditions which may insure the existence of an infinite
independent set  or consequences of the inexistence
 of such sets  have been considered within the framework of the structure of closure systems \index{closure system}
 (cf. the research on the  "free-subset problem" of Hajnal
\cite {shel} or on the cofinality of posets \cite{galv, mp2}). A
basic result is the following. \\
\begin{theorem} \cite {chapou}  \cite {lmp3} \label  {tm2.1} Let $\kappa$ be a cardinal number; for a join-semilattice $P$
the following properties are equivalent:\\ $(i)$ $P$ contains an
independent set of size $\kappa$;\\ $(ii)$ $P$ contains  a
join-subsemilattice isomorphic to $[\kappa]^{<\omega}$;\\ $(iii)$
$P$ contains  a subposet isomorphic to $[\kappa]^{<\omega}$;\\
$(iv)$ $J(P)$ contains a subposet isomorphic to $\mathfrak P
(\kappa)$;\\ $(v)$ $\mathfrak P (\kappa)$ embeds in $J(P)$ via a
map preserving arbitrary joins.
\end{theorem}
Let $L(\alpha):= 1+(1\oplus J(\alpha))+1$ be the lattice made of
the direct sum \index{direct sum} of the one-element chain $1$ and
the chain $J(\alpha)$, ($\alpha$  finite or equal to $\omega^*$),
with top and bottom added.
\begin{figure}
\begin{center}
\includegraphics[width=1in]{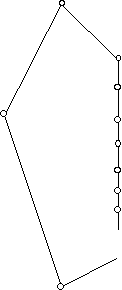}
\caption{$L(\omega^*)$}
\end{center}
\label{lomega}
\end{figure}

Clearly $J(\underline\Omega(\omega^*))$ contains a sublattice
isomorphic to  $L(\omega^*)$. Since a modular lattice
\index{modular lattice} contains no sublattice isomorphic  to
$L(2)$, we get as a corollary of Theorem \ref {thm4}:

\begin{theorem}\label{cor4}
An algebraic modular lattice $L$  is well-founded
\index{well-founded} if and only if $K(L)$ is well-founded  and
contains no infinite independent set.
\end{theorem}

 Another consequence is this:
\begin{theorem} \label{thmwf}
 For a join-semilattice $P$, the following properties are equivalent:
 \begin{enumerate}[(i)]
 \item $P$ is well-founded with no infinite antichain \index{antichain};
 \item $P$ contains no infinite independent set and embeds as a join-semilattice into a join-semilattice
 of the  form $I_{<\omega} (Q)$ where $Q$ is some well-founded poset.
\end{enumerate}
\end{theorem}

Posets which are well-founded and have no  infinite antichain are
said {\it well-partially-ordered}\index{well-partially-ordered} or
{\it well-quasi-ordered}\index{well-quasi-ordered}, wqo for short.
They play an important role in several areas (see \cite
{fraissetr}). If $P$ is a wqo join-semilattice then $J(P)$, the
poset of ideals  of  $P$,  is well-founded and one may assign to
every $J\in J(P)$ an ordinal, its  {\it height} \index{height},
denoted by $h(J, J(P))$. This ordinal is defined  by induction,
setting $h(J,J(P)):= Sup(\{h(J',J(P))+1: J'\in J(P), J'\subset
J\})$ and $h(J', J(P)):= 0$ if $J'$ is minimal in $J(P)$. The
ordinal $h(J(P)):= h(P, J(P))+1$ is the {\it height} of $J(P)$. If
$P:= I_{<\omega}(Q)$, with $Q$ wqo, then $J(P)$ contains a chain
of order type $h(J(P))$. This is an equivalent form of the famous
result of de Jongh and Parikh \cite{dejongh-parikh} asserting that
among the linear extensions \index{linear extension} of a wqo, one
has a maximum order type.

\begin{problem} Let $P$ be a wqo join-semilattice; does $J(P)$ contain a chain of order type $h(J(P))$?\end{problem}

An immediate corollary of Theorem \ref {thmwf} is:

\begin{corollary} \label{poset} A join-semilattice $P$ of $[\omega]^{<\omega}$ contains either  $[\omega]^{<\omega}$ as a join-semilattice or is wqo.
\end{corollary}
 Let us compare  join-subsemilattices of $[\omega]^{<\omega}$.
 Set $P\leq P'$ for two  such join-subsemilattices  if $P$ embeds in $P'$ as a join-semilattice.
 This gives a quasi-order  and,  according to  Corollary \ref{poset},
 the poset corresponding to this quasi-order has a largest element (namely $[\omega]^{<\omega}$),
 and  all other members come from   wqo join-semilattices.
 Basic examples of join-subsemilattices of $[\omega]^{<\omega}$ are the  $I_{<\omega} (Q)$'s
 where $Q$ is a countable poset such that no element is above infinitely many
 elements. These posets $Q$ are exactly those which are embeddable
 in the poset $[\omega]^{<\omega}$ ordered by inclusion. An
 interesting subclass is made of posets of the form $Q=(\N, \leq)$
 where the order $\leq$ is the intersection of the natural order $\mathfrak{N}$ on $\N$
 and of a linear order\index{linear order} $\mathfrak{L}$ on $\N$,
 (that is $x \leq y$ if $x \leq y$ w.r.t. $\mathfrak{N}$ and $x \leq y$ w.r.t. $\mathfrak{L}$).
 If $\alpha$ is the type of the linear order,
 a poset of this form  is a {\it sierpinskisation} \index{sierpinskisation} of $\alpha$.
 The corresponding join-semilattices are wqo provided that the posets  $Q$ have no infinite antichain;
 in the particular case of  a sierpinskisation of $\alpha$ this amounts to the fact that $\alpha$ is
 well-ordered \index{well-order}.

 As shown in \cite {pz}, sierpinskisations
given by a bijective map $\psi:\omega\alpha\rightarrow \omega$
which is  order-preserving \index{order-preserving} on each
component $\omega \cdot \{i\}$ of $\omega\alpha$ are all
embeddable in each other, and for this reason  denoted by the same
symbol $\Omega(\alpha)$. Among the representatives of $\Omega
(\alpha)$, some are join-semilattices, and among them, join-subsemilattices
of the direct product \index{direct product} $\omega\times \alpha$
(this is notably the case of the poset $\Omega(\omega ^{*})$ we
previously defined). We extend the first part of Proposition
\ref{w-f}, showing that except for $\alpha\leq\omega$, the
representatives of $\Omega(\alpha)$ which are join-semilattices
never embed in $[\omega]^{<\omega}$ as join-semilattices, whereas
they embed  as posets (see Corollary \ref{sierpinski} and Example
\ref{ex:ordinal}). From this result, it follows that the posets
$\Omega (\alpha)$ and $I_{<\omega}(\Omega (\alpha))$  do not embed
in each other as join-semilattices.

 These two posets  provide  examples of a  join-semilattice
$P$ such that $P$ contains no chain of type $\alpha$ while $J(P)$
contains a chain of type $J(\alpha)$. However,   if $\alpha$ is
not well ordered  then $I_{<\omega}(\Omega (\alpha))$ and
$[\omega]^{<\omega}$ embed in each other as join-semilattices.

  \begin{problem}
 Let $\alpha$ be a countable  ordinal.
Is there a minimum member among the join-subsemilattices $P$ of
$[\omega]^{<\omega}$ such that $J(P)$ contains a chain of type
$\alpha+1$? Is it true that this minimum is
$I_{<\omega}(\Omega(\alpha))$  if $\alpha$ is indecomposable?
\end{problem}

\section{Definitions and basic results}\label {subsection2.2}
Our definitions and notations are standard and agree with \cite
{grat} except on minor points that we will mention. We adopt the
same terminology as in \cite {chapou}. We recall only few things.
Let $P$ be a poset. A subset $I$ of $P$ is an {\it initial
segment} \index{initial segment} of $P$ if $x\in P$, $y\in I$ and
$x\leq y$ imply $x\in I$. If
 $A$ is a subset of $P$, then $\downarrow A=\{x\in P: x\leq y$ for some $y\in A\}$ denotes the least initial segment containing
$A$. If $I=\downarrow A$ we say that $I$ is {\it generated} by $A$
or $A$ is {\it cofinal} in \index{cofinal} $I$. If $A=\{a\}$ then
$I$ is a {\it principal initial segment} \index{principal initial
segment} and we write $\downarrow a$ instead of $\downarrow
\{a\}$. We denote $down(P)$ the set of principal initial segments
of $P$. A {\it final segment} \index{final segment} of $P$ is any
initial segment of $P^*$, the dual of $P$. We denote by $\uparrow
A$ the final segment generated by $A$. If $A=\{a\}$ we write
$\uparrow a$ instead of $\uparrow \{a\}$. A subset $I$ of $P$ is
{\it directed}\index{directed}
 if every finite subset of $I$
has an upper bound in $I$ (that is $I$ is non-empty and every pair
of elements of $I$ has an upper bound). An {\it ideal}\index{ideal}
is a non-empty directed initial segment of $P$ (in some other texts,
the empty set is an ideal). We denote $I(P)$ (respectively, $I_
{<\omega }(P)$, $J(P)$) the set of initial segments (respectively,
finitely generated initial segments, ideals of $P$) ordered by
inclusion and we set
 $ J_ *(P):=J(P)\cup\{\emptyset\}$, $I_ 0 (P):=I_ {<\omega}
(P)\setminus\{\emptyset\}$. Others authors use {\it down set} for
initial segment. Note that $down(P)$ has not to be confused with
$I(P)$. If $P$ is a join-semilattice with a $0$, an element $x\in
P$ is {\it join-irreducible} \index{join-irreducible} if it is
distinct from $0$, and if $x=a\vee b$ implies $x=a$ or $x=b$ (this
is a slight variation from \cite{grat}). We denote $\J_{irr}(P)$
the set of join-irreducibles of $P$. An element $a$ in a lattice
$L$ is {\it compact}\index{compact} if for every $A\subset L$,
$a\leq \bigvee A$ implies $a\leq \bigvee A'$ for some finite
subset $A'$ of $A$. The lattice $L$ is {\it compactly
generated}\index{ compactly generated} if every element is a
supremum of compact elements. A lattice is {\it
algebraic}\index{algebraic} if it is complete and compactly
generated.
\\

We note that $I_{ <\omega}(P)$ is the set of compact elements of
$I(P)$, hence $J(I_{ <\omega}(P))\cong I(P)$. Moreover $I_
{<\omega}(P)$ is  a lattice, and in fact a distributive lattice,
if and only if $P$ is {\it $\downarrow$-closed}
\index{$\downarrow$-closed}, that is, the intersection of two
principal initial segments of $P$ is a finite union, possibly
empty, of principal initial segments. We also note that $J(P)$  is
the set of join-irreducible elements of $I(P)$; moreover,
$I_{<\omega}(J(P))\cong I(P)$ whenever $P$ has no infinite
antichain. \\ Notably for the proof of Theorem
\ref{finitesubsets}, we will need the following  results.

\begin{theorem} \label{wellfounded} Let $P$ be a poset.\\
$a)$  $I_{ <\omega}(P)$ is well-founded if and only if $P$ is
well-founded (Birkhoff  1937, see \cite {birk});\\ $b)$ $I_{
<\omega}(P)$ is wqo iff $P$ is wqo iff $I(P)$ well-founded (
Higman 1952 \cite {higm}); \\ $c)$ if $P$ is a well-founded
join-semilattice with a 0, then every member of $P$ is a finite
join of join-irreducible elements of $P$ (Birkhoff, 1937, see
\cite {birk});\\ $d)$ A join-semilattice $P$ with a zero is wqo if
and only if every member of $P$ is a finite join of
join-irreducible elements of $P$ and the set $\J_{irr}(P)$ of
these join-irreducible elements is wqo (follows from $b)$ and
$c)$).
\end{theorem}
 A poset $P$ is {\it scattered} \index{scattered} if it does not contain  a copy of $\eta$,
 the chain of rational numbers. A topological space $T$ is
{\it scattered} if every non-empty closed set contains some
isolated point. The power set of a set, once  equipped with the
product topology, is a compact space. The set $J(P)$ of ideals of
a join-semilattice $P$ with a $0$ is a closed subspace of
$\mathfrak P (P)$, hence is a compact space too. Consequently, an
algebraic lattice $L$ can be viewed as a poset and a topological
space as well. It is easy  to see that if $L$ is topologically
scattered \index{topologically scattered} then it is order
scattered \index{order scattered}. It is  a more significant fact,
due to M.Mislove \cite {misl}, that the converse holds if $L$ is
distributive.

\section[Separating chains of ideals]
{Separating chains of ideals and proofs of Proposition \ref {w-f}
and Theorem \ref{thm4}}

Let $P$ be a join-semilattice. If $x\in P$ and $J\in J(P)$, then
$\downarrow x$ and $J$ have a join $\downarrow x \bigvee J$ in
$J(P)$ and $\downarrow x \bigvee J=\downarrow\{ x \vee y: y\in J\}$.
Instead of $\downarrow x \bigvee J$ we also use the notation
$\{x\}\bigvee J$. Note that $\{x\}\bigvee J$ is the least member of
$J(P)$ containing $\{x\}\cup J$. We say that a non-empty chain
$\mathcal I$ of ideals of $P$ is {\it separating} \index{separating}
if for every $I\in\mathcal I\setminus\{\cup\mathcal I\}$ and every
$x\in \cup\mathcal I\setminus I$, there is some $J\in\mathcal I$
such that $I\not\subseteq\{x\}\bigvee J$.\\ If $\mathcal I$ is
separating then $\mathcal I$ has a least element implies it is a
singleton set. In $P:=[\omega]^{<\omega}$, the chain $\mathcal
I:=\{I_{n}: n<\omega\}$ where $I_{n}$ consists of the finite subsets
of $\{m: n\leq m\}$ is separating. In $P:=\omega^*$, the chain
$\mathcal I:=\{\downarrow x: x\in P\}$ is non-separating, as well as
all of its infinite subchains. In $P:=\Omega(\omega^*)$ the chain
$\mathcal I:=\{I_{n}: n<\omega\}$ where $I_{n}:=\{(i,j): n\leq
i<j<\omega\}$ has the same property.\\ We may observe that {\it a
join-preserving embedding\index{join-preserving embedding} from a
join-semilattice $P$ into a join-semilattice $Q$ transforms every
separating (resp. non-separating) chain of ideals of $P$ into a
separating (resp. non-separating) chain of ideals of $Q$} (If
$\mathcal I$ is a separating chain of ideals of $P$, then $\mathcal
J=\{f(I): I\in \mathcal I\}$ is a separating chain of ideals of
$Q$). Hence the containment of $[\omega]^{<\omega}$ (resp. of
$\omega^*$ or of $\Omega(\omega^*)$), as a join-subsemilattice,
provides a chain of ideals which is separating (resp.
non-separating, as are all its infinite subchains, as well). We show
in the  next  two lemmas that the converse holds.

\begin{lemma}\label{independent}
A join-semilattice $P$ contains an infinite independent set if and
only if it contains an infinite separating chain of ideals.
\end{lemma}

\begin{proof}
Let $X=\{x_{n}: n<\omega\}$ be an infinite independent set. Let
$I_{n}$ be the ideal generated by $X\setminus\{x_{i}: 0\leq i\leq
n\}$. The chain $\mathcal I=\{I_{n}: n<\omega\}$ is separating. Let
$\mathcal I$ be an infinite separating chain of ideals. Define
inductively an infinite sequence $x_{0}, I_{0}, \ldots, x_{n},I_{n},
\ldots$ such that $I_{0}\in\mathcal I\setminus\{\cup\mathcal I\},
x_{0}\in \cup\mathcal I \setminus I_{0}$ and
 such that:\\
 $a_{n})$
$I_{n}\in\mathcal I$;\\ $b_{n})$ $I_{n}\subset I_{n-1}$;\\
 $c_{n})$
$x_{n}\in I_{n-1}\setminus(\{x_{0}\vee \ldots \vee
x_{n-1}\}\bigvee I_{n})$ for every  $n\geq 1$.\\ The construction
is immediate. Indeed, since $\mathcal I$ is infinite then
$\mathcal I\setminus\{\cup\mathcal I\}\not =\emptyset$. Choose
arbitrary $I_{0}\in\mathcal I\setminus\{\cup\mathcal I\}$ and
$x_{0}\in\cup\mathcal I\setminus I_{0}$. Let $n\geq 1$. Suppose
$x_{k}, I_{k}$ defined and satisfying $a_{k}), b_{k}), c_{k})$ for
all $k\leq n-1$. Set $I:=I_{n-1}$ and $x:=x_{0}\vee \ldots \vee
x_{n-1}$. Since $I\in\mathcal I$ and $x\in\cup\mathcal I\setminus
I$, there is some $J\in\mathcal I$ such that
$I\not\subseteq\{x\}\bigvee J$. Let $z\in I\setminus(\{x\}\bigvee
J)$. Set $x_{n}:=z$, $I_{n}:=J$. The set $X:=\{x_{n}: n<\omega\}$
is independent. Indeed if $x\in X$ then since $x=x_{n}$ for some
$n$, $n<\omega$, condition $c_{n})$ asserts that there is some
ideal containing $X\setminus \{x\}$ and excluding $x$.
\end{proof}\\

\begin{lemma} \label{w*}
A join-semilattice $P$ contains either $\omega^*$ or
$\Omega(\omega^*)$ as a join-subsemilattice if and only if it
contains an $\omega^*$-chain $\mathcal I$ of ideals such that all
infinite subchains are non-separating.
\end{lemma}

\begin{proof}
Let $\mathcal I$ be an $\omega^*$-chain of ideals and let $A$ be its
largest element (that is $A=\cup\mathcal I$). Let $E$ denote the set
$\{x: x\in A$ and $I\subset\downarrow x$ for some $I\in \mathcal
I\}$.\\ {\bf Case (i)}. For every $I\in\mathcal I$, $I\cap
E\not=\emptyset$. We can build an infinite strictly decreasing
sequence $x_{0}, \ldots, x_{n}, \ldots$ of elements of $P$. Indeed,
let us choose $x_{0}\in E\cap(\cup\mathcal I)$ and $I_{0}$ such that
$I_{{0}}\subset\downarrow x_{0}$. Suppose $x_{0}, \ldots, x_{n}$ and
$I_{0}, \ldots, I_{n}$ defined such that $I_{i}\subset\downarrow
x_{i}$ for all $i=0, \ldots, n$. As $E\cap I_{n}\not=\emptyset$ we
can select $x_{n}\in E\cap I_{n}$ and by definition of $E$, we can
select some $I_{n+1}\in \mathcal I$ such
that $I_{n+1}\subset\downarrow x_{n+1}$. Thus $\omega^*\leq P$.\\
{\bf Case (ii)}. There is some $I\in\mathcal I$ such that $I\cap
E=\emptyset$. In particular all members of $\mathcal I$ included in
$I$ are unbounded in $I$. Since all infinite subchains of $\mathcal
I$ are non-separating then, with no loss of generality, we may
suppose that $I=A$ (hence $E=\emptyset$). We set $I_{-1}:=A$ and
define a sequence $x_{0}, I_{0}, \ldots, x_{n}, I_{n}, \ldots$ such
that
 $I_{n}\in\mathcal I$, $x_{n}\in I_{n-1}\setminus I_{n}$ and $I_{n}\subseteq \{x_{n}\}\bigvee I$ for all $I\in\mathcal I$, all $n<\omega$.
Members of this sequence being defined for all $n', n'<n$, observe
that the set $\mathcal I_{n}:=\{I\in\mathcal I: I\subseteq
I_{n-1}\}$ being infinite is non-separating, hence there are $I\in
\mathcal I_{n}$ and $x\in I_{n-1}\setminus I$ such that $I\subseteq
\{x\}\bigvee J$ for all $J\in\mathcal I_{n}$. Set $I_{n}:=I$ and
$x_{n}:=x$. Next, we define a sequence $y_{0}:=x_{0}, \ldots, y_{n},
\ldots$ such that for every $n\geq 1$:\\ $a_{n})$ $x_{n} \leq
y_{n}\in I_{n-1}$; \\$b_{n})$ $y_{n} \not\leq y_{0}\vee y_{n-1}$;\\
$c_{n})$  $y_{j}\leq y_{i}\vee y_{n}$ for every $i\leq j\leq n$.\\
Suppose $y_{0}, \ldots, y_{n-1}$ defined for some $n$, $n\geq 1$.
Since $I_{n-1}$ is unbounded, we may select $z\in I_{n-1}$ such that
$z\not\leq y_{0}\vee\ldots\vee y_{n-1}$. If $n=1$, we set
$y_{1}:=x_{1}\vee z$. Suppose $n\geq 2$. Let $0\leq j\leq n-2$.
Since $y_{j+1}\vee\ldots\vee y_{n-1}\in I_{j}\subseteq \{
x_{j}\}\bigvee I_{n-1}$  we may select $t_{j}\in I_{n-1}$ such that
$y_{j+1}\vee\ldots\vee y_{n-1}\leq x_{j}\vee t_{j}$. Set
$t:=t_{0}\vee\ldots\vee t_{n-2}$ and $y_{n}:=x_{n}\vee z\vee t$.\\
Let $f: \Omega(\omega^*)\to P$ be defined by
$f(i,j):=y_{i}\vee y_{j}$ for all $(i,j)$, $i<j<\omega$.\\
Condition $c_{n})$ insures that $f$ is join-preserving. Indeed, let
$(i,j), (i',j')\in\Omega(\omega^*)$. We have $(i,j)\vee
(i',j')=(i\wedge i',j\vee j')$ hence $f((i,j)\vee (i',j'))=f(i\wedge
i',j\vee j')=y_{i\wedge i'}\vee y_{j\vee j'}$. If $F$ is a finite
subset of $\omega$ with minimum $a$ and maximum $b$ then conditions
$c_{n})$ force $\bigvee\{y_{n}: n\in F\}=y_{a}\vee y_{b}$. If
$F:=\{i,j,i',j'\}$ then, taking account of $i<j$ and $i'<j'$, we
have $f(i,j)\vee f(i',j')=y_{i}\vee y_{j}\vee y_{i'}\vee
y_{j'}=y_{i\wedge i'}\vee y_{j\vee j'}$. Hence
$f((i,j)\vee (i',j'))=f(i,j)\vee f(i',j')$, proving our claim.  \\
Next, $f$ is one-to-one. Let $(i,j), (i',j')\in\Omega(\omega^*)$
such that $f(i,j)=f(i',j')$, that is $y_{i}\vee y_{j}=y_{i'}\vee
y_{j'}$ $(1)$. Suppose $j<j'$. Since $0\leq i<j$, Condition $c_{j})$
implies $y_{i}\leq y_{0}\vee y_{j}$. In the other hand,  since
$0\leq j\leq j'-1$, Condition $c_{j'-1})$ implies $y_{j}\leq
y_{0}\vee y_{j'-1}$. Hence $y_{i}\vee y_{j}\leq y_{0}\vee y_{j'-1}$.
From $(1)$ we get $y_{j'}\leq y_{0}\vee y_{j'-1}$, contradicting
Condition $b_{j'})$. Hence $j'\leq j$. Exchanging the roles of
$j,j'$ gives $j'\leq j$ thus $j=j'$. If $i<i'$ then, Conditions
$a_{i'})$ and $a_{j'})$ assure $y_{i'}\in I_{i'-1}$ and $y_{j'}\in
I_{j'-1}$. Since $I_{j'-1}\subseteq I_{i'-1}$ we have $y_{i'}\vee
y_{j'}\in I_{i'-1}$. In the other hand $x_{i}\not\in I_{i}$ and
$x_{i}\leq y_{i}\vee y_{j}$ thus $y_{i}\vee y_{j}\not\in I_{i}$.
From $I_{i'-1}\subseteq I_{i}$, we have $y_{i}\vee y_{j}\not\in
I_{i'-1}$, hence $y_{i}\vee y_{j}\not = y_{i'}\vee y_{j'}$ and
$i'\leq i$. Similarly we get $i\leq i'$. Consequently $i=i'$.
\end{proof}

\subsection{Proof of Proposition \ref{w-f}}
If $\underline \Omega(\omega^*)$ embeds in $[\omega]^{<\omega}$ then
$[\omega]^{<\omega}$ contains a non-separating $\omega^*$-chain of
ideals. This is impossible: a  non-separating chain of ideals of
$[\omega]^{<\omega}$ has necessarily a least element. Indeed, if the
pair $x,I$ ($x\in [\omega]^{<\omega}$, $I\in  \mathcal I$) witnesses
the fact that the chain $\mathcal I$ is non-separating  then there
are at most $\mid x\mid+1$ ideals belonging to  $\mathcal I$ which
are included in $I$ (note that the set $\{\cup I \setminus \cup J:
J\subseteq I, J\in \mathcal I\}$  is a chain of subsets of $x$). The
proof of the general case requires more care. If $\underline
\Omega(\omega^*)$ embeds in $I_{<\omega}(Q)$ as a join-semilattice
then we may find a sequence $x_{0}, I_{0}, \ldots, x_{n}, I_{n},
\ldots$ such that $I_{n}\subset I_{n-1}\in J(I_{<\omega}(Q))$,
$x_{n}\in I_{n-1} \setminus I_{n}$ and $I_{n} \subseteq
\{x_{n}\}\bigvee I_{m}$ for every $n<\omega$ and every $m<\omega$.
Set $I_{\omega}:=\bigcap\{ I_{n}: n<\omega\}$, $\overline
I_{n}:=\cup I_{n}$ for every $n\leq\omega$, $Q':=Q\setminus
\overline I_{\omega}$ and $y_{n}:=x_{n}\setminus \overline
I_{\omega}$ for every $n<\omega$. We claim that $y_{0}, \ldots,
y_{n}, \ldots$ form a strictly descending sequence in
$I_{<\omega}(Q')$. According to Property $a)$  stated in Theorem
\ref {wellfounded},  $Q'$, thus $Q$, is not well-founded.\\ First,
$y_{n}\in I_{<\omega}(Q')$. Indeed, if $a_{n}\in [Q]^{<\omega}$
generates $x_{n}\in I_{<\omega}(Q)$ then, since $\overline
I_{\omega}\in I(Q)$, $a_{n}\setminus \overline I_{\omega}$ generates
$x_{n}\setminus \overline I_{\omega}\in I(Q')$. Next,
$y_{n+1}\subset y_{n}$. It suffices to prove that the following
inclusions hold: $$x_{n+1}\cup \overline I_{\omega}\subseteq
\overline I_{n}\subset x_{n}\cup \overline I_{\omega}$$ Indeed,
substracting $\overline I_{\omega}$, from the sets  figuring above,
we get:\\ $$y_{n+1}=(x_{n+1}\cup \overline I_{\omega})\setminus
 \overline I_{\omega}\subset (x_{n}\cup \overline I_{\omega})\setminus \overline I_{\omega}=y_{n}$$
The first inclusion is obvious. For the second note that, since
$J(I_{<\omega}(Q))$ is isomorphic to $I(Q)$, complete distributivity
holds, hence  with the hypotheses on the sequence $x_{0}, I_{0},
\ldots, x_{n}, I_{n}, \ldots$  we have $I_{n} \subseteq \bigcap\{\{
x_{n}\}\bigvee I_{m}: m<\omega\}=\{ x_{n}\}\bigvee\bigcap\{ I_{m}:
m<\omega\}=\{ x_{n}\}\bigvee I_{\omega}$, thus $\overline
I_{n}\subset x_{n}\cup \overline I_{\omega}$.
\endproof\\
\begin{remark} One can deduce the fact that $\Omega(\omega^*)$ does not embed as  a join-semilattice in
$[\omega]^{<\omega}$ from the fact that it contains a strictly
descending chain of completely meet-irreducible \index{completely
meet-irreducible} ideals (namely the chain $\mathcal I:=\{I_{n}:
n<\omega\}$ where $I_{n}:=\{(i,j): n\leq i<j<\omega\}$) (see
Proposition \ref {lem03}) but this fact by itself does not prevent
the existence of some well-founded poset $Q$ such that
$\Omega(\omega^*)$ embeds as a join semilattice in
$I_{<\omega}(Q)$. \end{remark}

\subsection{Proof of Theorem \ref{thm4}}
In terms of join-semilattices and  ideals, result becomes this:
let $P$ be a join-semilattice, then $J(P)$ is well-founded if and
only if $P$ is well-founded  and contains no join-subsemilattice
isomorphic to $\Omega(\omega^*)$ or to $[\omega]^{<\omega}$.\\

The proof goes as follows. Suppose that $J(P)$ is not well-founded. If some $\omega^*$-chain in $J(P)$ is separating then, according to
Lemma \ref{independent}, $P$
contains an infinite independent set. From Theorem \ref {tm2.1}, it contains a join-subsemilattice isomorphic to $[\omega]^{<\omega}$.
If no $\omega^*$-chain in $J(P)$ is separating, then all the infinite subchains of an arbitrary $\omega^*$-chain are non-separating.
From Lemma \ref{w*},
 either $\omega^*$ or
$\Omega ( \omega^*)$ embed in $P$ as a join-semilattice. The
converse is obvious. \endproof\\

\section[Join-subsemilattices of $I_{<\omega}(Q)$]{Join-subsemilattices of $I_{<\omega}(Q)$ and proof of Theorem  \ref{thmwf}}
In this section, we consider  join-semilattices  which embed in
join-semilattices of the form $I_{<\omega}(Q)$. These are easy  to
characterize internally (see Proposition \ref {lem04}). This is
also the case if the posets $Q$ are antichains (see Proposition
\ref {lem03}) but does not go so well if the posets  $Q$ are
well-founded (see Lemma \ref{counterexample}).

Let us recall  that if $P$ is a  join-semilattice, an element $x\in
P$ is  {\it join-prime}\index{join-prime} (or prime if there is no
confusion), if it is distinct from the least element $0$, if any,
and if $x\leq a\vee b$ implies $x\leq a$ or $x\leq b$. This amounts
to the fact that $P\setminus \uparrow x$ is an ideal . We denote
$\J_{pri}(P)$, the set of join-prime  members of $P$. We recall that
$\J_{pri}(P)\subseteq \J_{irr}(P)$; the equality holds provided that
$P$ is a distributive lattice. It also holds if $P=
{I}_{<\omega}(Q)$. Indeed:

\begin{fact}For  an arbitrary poset $Q$, we have:
\begin{equation}\label{eqirr}
 {\J}_{irr}( I_{<\omega} (Q))=\J_{pri}( I_{<\omega} (Q))= down (Q)
 \end{equation}
\end{fact}
\begin{fact} For a poset $P$, the following properties are equivalent:
\begin{itemize}
\item $P$ is isomorphic to $ {I}_{<\omega}(Q)$ for some poset $Q$;
\item $P$ is a join-semilattice with a least element  in which  every element is a finite  join of primes.
\end{itemize}
\end{fact}

\begin{proof} Observe that the primes in $I_{<\omega} (Q)$, are the $\downarrow x$, $x\in
Q$. Let $I\in I_{<\omega} (Q)$ and $F\in [Q]^{<\omega}$ generating
$I$, we have $I=\cup \{\downarrow x: x\in F\}$ . Conversely, let
$P$ be a join-semilattice with a $0$. If every element in $P$ is a
finite join of primes, then $P\cong I_{<\omega} (Q)$ where
$Q:=\J_{pri}(P)$.
\end{proof}

Let $L$ be  a complete lattice \index{complete lattice}. For $x\in
L$ , set $x^{+}:=\bigwedge\{y\in L: x<y\}$. We recall that $x\in L$
is {\it completely meet-irreducible} \index{completely
meet-irreducible} if $x=\bigwedge X$ implies $x\in X$, or
-equivalently- $x\neq x^{+}$.  We denote $\triangle (L)$ the set of
completely meet-irreducible members of $L$. We recall the following
Lemma.
\begin{lemma}\label{lem01} Let $P$ be a join-semilattice,  $I\in J(P)$ and $x\in P$.
Then $x\in I^{+}\setminus I$ if and only if  $I$ is a maximal ideal of $P\setminus
\uparrow x$.
\end{lemma}

\begin{proposition}\label{lem04} Let $P$ be a join-semilattice. The
following properties are equivalent:
\begin {enumerate}[(i)]
\item $P$ embeds in
$I_{<\omega}(Q)$,  as a join-semilattice, for some poset $Q$;
\item $P$ embeds in $I_{<\omega}(J(P))$ as a join-semilattice;
\item $P$ embeds in $I_{<\omega}(\triangle(J(P)))$ as a join-semilattice;
\item  For every $x\in P$, $P\setminus \uparrow x$ is a finite
union of ideals.
\end{enumerate}
\end{proposition}

\begin{proof}
$(i)\Rightarrow (iv)$ Let $\varphi$ be an embedding from $P$ in
$P':= I_{<\omega}(Q)$. We may suppose that $P$ has a least element
$0$ and that $\varphi (0)=\emptyset$ (if $P$ has no least element,
add one, say $0$,  and set $\varphi (0):=\emptyset $; if $P$ has a
least element, say $a$,  and $\varphi (a)\not =\emptyset$,  add to
$P$ an element $0$ below $a$ and set $\varphi (0):= \emptyset$).
For $J'\in \mathfrak { P}(P')$,  let $\varphi^{-1}(J'):= \{x\in P:
\varphi (x)\in J'\}$. Since $\varphi$ is order-preserving,
$\varphi^{-1}(J')\in I(P)$ whenever $J'\in I(P')$ ; moreover,
since $\varphi$ is join-preserving, $\varphi^{-1}(J')\in J(P)$
whenever $J'\in J(P')$. Now, let  $x\in P$. We have
$\varphi^{-1}(P'\setminus \varphi(x)):= P\setminus \uparrow x$.
Since $\varphi(x)$ is a finite join of primes, $P'\setminus
\uparrow  \varphi (x)$ is a finite union of ideals. Since their
inverse images are ideals, $P\setminus \uparrow x$ is a finite
union of ideals too.

$(iv)\Rightarrow (iii)$  We use  the well-known method for representing a poset by a family of sets.

\begin{fact}\label{representation}Let $P$ be a poset and $Q\subseteq I(P)$. For $x\in P$ set
$\varphi _{Q}(x):=\{J\in Q: x\not \in J\}$.  Then:
\begin{enumerate}[(a)]
\item  $\varphi _{Q}(x) \in I(Q)$;
\item  $\varphi _{Q}:P\rightarrow I(Q)$ is an order-preserving \index{order-preserving} map;
\item $\varphi _{Q}$  is  an order-embedding\index{order-embedding} if and only if for every $x, y\in P$ such that $x\not \leq y$
there is some $J\in Q$ such that $x\not \in J$ and $y\in J$.
\end{enumerate}
\end{fact}

Applying  this to $Q:=\triangle( J(P))$ we get  immediately that
$\varphi _{Q}$ is join-preserving \index{join-preserving}. Moreover,
$\varphi _{Q}(x)\in I_{<\omega}(Q)$ if and only if $P\setminus
\uparrow x$ is a finite union of ideals. Indeed, we have
$P\setminus \uparrow x= \cup \varphi _{Q}(x)$, proving that
$P\setminus \uparrow x$ is a finite union of ideals provided that
$\varphi _{Q}(x)\in I_{<\omega}(Q)$. Conversely, if $P\setminus
\uparrow x$ is a finite union of ideals, say $I_{0}, \dots, I_{n}$,
then since ideals are prime members of $I(P)$, every ideal included
in $I$ is included in some $I_i$, proving that $\varphi _{Q}(x)\in
I_{<\omega}(Q)$. To conclude, note that if $P$ is a join-semilattice
then $\varphi _{Q}$ is join-preserving.

$(iii)\Rightarrow (ii)$ Trivial.

$(ii)\Rightarrow (i)$ Trivial.

\end{proof}

\begin{corollary}\label{corthm16} If a   join-semilattice $P$ has no infinite antichain, it embeds in
$I_{<\omega}(J(P))$ as a join-subsemilattice.
\end{corollary}
\begin{proof} As is well known, if  a poset has no infinite antichain then every initial segment
is a finite union of ideals (cf \cite{Erdos-Tarski}, see also
\cite{fraissetr} 4.7.3 pp. 125). Thus Proposition \ref{lem04}
applies.
\end{proof}

Another corollary of Proposition  \ref{lem04}  is the following.

\begin{corollary}\label{lem00}  Let $P$ be a join-semilattice. If  for every $x\in P$,  $P\setminus \uparrow x$ is a finite union of ideals and $ \triangle (J(P))$ is well-founded   then   $P$ embeds  as a join-subsemilattice in $I_{<\omega}(Q)$,  for some well-founded poset $Q$.
\end{corollary}

The converse does not hold:

\begin{example}\label{counterexample} There is a bipartite poset $Q$ such that $I_{<\omega}(Q)$ contains a
join-semilattice $P$ for which  $ \triangle (J(P))$ is not
well-founded.
\end{example}

\begin{proof}
Let $\underline 2:= \{0,1\}$ and $Q:=\mathbb{N}\times
\underline{2}$. Order $Q$ in such a way that $(m,i)<(n,j)$ if
$m>n$ in $\N$ and $i<j$ in $\underline 2$.

Let $P$ be the set
of subsets $X$ of $Q$ of the form $X:=F\times \{0\}\cup G\times
\{1\}$ such that $F$ is a  non-empty final segment of $\mathbb{N}$, $G$ is a
non-empty finite subset of $\mathbb{N}$ and
\begin{equation}\label{contrex}
min(F)-1\leq min(G)\leq
min(F)
\end{equation}
where $min(F)$ and $min(G)$  denote the least element of $F$  and
$G$ w.r.t. the natural order on $\mathbb{N}$. For each $n\in \N$,
let $I_{n}:=\{ X\in P:  (n,0) \not \in X\}$.

{\bf Claim }
\begin{enumerate}
\item $Q$ is bipartite and  $P$ is a join-subsemilattice of $I_{<\omega}(Q)$.
\item The $I_{n}$'s form a strictly descending sequence \index{descending sequence} of members of
$ \triangle (J(P))$.
\end{enumerate}

{\bf Proof of the Claim }

1. The poset $Q$ is decomposed into two antichains, namely $\N\times
\{0\}$ and $\N\times \{1\}$ and for this raison is called {\it
bipartite}.
 Next, $P$ is a
subset of $I_{<\omega}(Q)$. Indeed, Let $X\in P$. Let $F, G$ such
that $ X=F\times \{0\}\cup G\times \{1\}$.  Set $G':=
G\times\{1\}$. If $min(G)=min(F)-1$, then $X=\downarrow G'$
whereas if  $min(G)=min(F)$ then $X=\downarrow G' \cup \{(min(F),
0)\}$. In both cases    $X\in I_{<\omega}(Q)$.  Finally, $P$ is a
join-semilattice. Indeed,  let $X,X'\in P$ with $X:=F\times
\{0\}\cup G\times \{1\}$ and $X':=F'\times \{0\}\cup G'\times
\{1\}$. Obviously $X\cup X'=(F\cup F')\times \{0\}\cup (G\cup
G')\times \{1\}$. Since $X, X'\in P$, $F\cup F'$ is a non-empty
final segment of $\mathbb{N}$ and $G\cup G'$ is a non-empty finite
subset of $\mathbb{N}$. We have $min(G\cup
G')=min(\{min(G),min(G')\})\leq min(\{min(F),min(F')\})=min(F\cup
F')$ and similarly $min(F\cup F')-1=min\{min(F),
min(F')\}-1=min\{min(F)-1, min(F')-1\} \leq min
\{min(G),min(G')\}=min(G\cup G')$, proving that inequalities as in
(\ref {contrex}) hold. Thus $X\cup X'\in I_{<\omega}(Q)$.

2.  Due to its definition, $I_{n}$ is an non-empty initial segment
of $P$ which is closed under finite unions, hence $I_{n}\in J(P)$.
Let $X_{n}:= \{(n,1), (m, 0): m\geq n+1\} $ and $Y_{n}:= X_{n}\cup
\{(n,0)\}$. Clearly, $X_{n} \in I_{n}$ and $Y_{n} \in P$. We claim
that $I_{n}^{+}=I_{n}\bigvee \{Y_{n}\}$. Indeed, let $J$ be an
ideal containing strictly $I_{n}$. Let $Y:=\{m\in \mathbb{N}:
m\geq p\}\times\{0\}\cup G\times\{1\} \in J\setminus I_{n}$. Since
$Y\not\in I_{n}$, we have  $p\leq n$ hence  $Y_{n}\subseteq Y\cup
X_{n}\in J$. It follows that $Y_{n}\in J$, thus
$I_{n}^{+}\subseteq J$, proving our claim. Since $I_{n}^{+}\not=
I_{n}$, $I_{n}\in\triangle (J(P))$. Since, trivially,
$I_{n}^{+}\subseteq I_{n-1}$ we have $I_{n}\subset I_{n-1}$,
proving that the $I_n$'s form a strictly descending sequence.
\end{proof}

Let  $E$ be a set and $\mathcal{F}$ be a subset of $\mathfrak
P(E)$, the power set of $E$. For $x\in E$, set $\mathcal {F}_{\neg
x}:=\{F\in \mathcal{F}: x\not\in F\}$ and for $X \subset \mathcal
F$ , set $\overline X:= \bigcup X$. Let $\mathcal{F}^{<\omega}$
(resp. $\mathcal{F}^\cup$) be the collection of finite (resp.
arbitrary) unions of members of $\mathcal{F}$. Ordered by
inclusion, $\mathcal{F}^\cup$ is a  complete lattice
\index{complete lattice}, the least element and the largest
element being  the empty set and $\bigcup\mathcal{F}$,
respectively.

\begin{lemma}\label{cl2} Let $Q$ be a poset,  $\mathcal{F}$ be  a subset of $I_{<\omega}(Q)$ and $P:= \mathcal{F}^{<\omega}$ ordered by inclusion.
\begin{enumerate}[{(a)}]

\item The map $X\rightarrow \overline X$ is an isomorphism  from  $J(P)$ onto  $\mathcal{F}^\cup$ ordered by inclusion.

\item  If $I\in \triangle (J(P))$
then there is
some $x\in Q$ such that $I=P_{\neg x}$.

\item If $\downarrow q$ is finite for every  $q\in Q$ then $\overline{I^+}\setminus\overline{I}$ is finite for every  $I\in J(P)$ and
the set $\varphi_{\triangle}(X):= \{I\in {\triangle}(J(P)) : X\not\in I\}$ is finite for every $ X\in P$.

\end{enumerate}
\end{lemma}

\begin{proof}

\noindent $(a)$ Let $I$ and $J$ be two ideals of $P$. Then $J$
contains $I$ if and only if $\overline{J}$ contains $\overline{I}$.
Indeed, if $I\subseteq J$ then, clearly $\overline{I}\subseteq
\overline{J}$. Conversely, suppose $\overline{I}\subseteq
\overline{J}$. If $X\in I$, then $X\subseteq \overline{I}$, thus
$X\subseteq \overline{J}$. Since $X\in I_{<\omega}(Q)$, and
$X\subseteq \overline{J}$, there are $X_{1}, \ldots, X_{n}\in J$
such that $X\subseteq Y=X_{1}\cup \ldots \cup X_{n}$. Since $J$ is
an ideal $Y\in J$. It follows that $X\in J$.

\noindent $(b)$ Let  $I\in \triangle (J(P))$. From $(a)$,  we have
$\overline{I}\subset\overline{I^+}$. Let $x\in
\overline{I^+}\setminus\overline{I}$.  Clearly  $P_{\neg x}$ is an
ideal containing $I$. Since  $x\not\in \overline{P_{\neg x}}$,
$P_{\neg x}$ is distinct from $I^{+}$. Hence $P_{\neg x}=I$. Note
that the converse of assertion $(b)$ does not  hold in general.

\noindent $(c)$ Let
$I\in \triangle (J(P))$ and $X\in I^{+}\setminus I$.  We have $\{X\}\bigvee
I= I^+$, hence  from$(a)$ $\overline{\{X\} \bigvee
I}= \overline{I^+}$. Since  $\overline{\{X\} \bigvee
I}=X\cup \overline{I}$ we have
$\overline{I^+}\setminus\overline{I}\subseteq X$. From our hypothesis on $P$,  $X$ is
finite, hence $\overline{I^+}\setminus\overline{I}$ is finite. Let $X\in P$.  If  $I\in \varphi_{\triangle}(X)$ then according  to $(b)$
there is some $x\in Q$ such that $I= P_{\neg x}$. Necessarily $x\in X$.
Since $X$ is finite, the number of these $I$'s is finite.
\end{proof}

\begin{proposition}\label{lem03} Let $P$ be a join-semilattice.
The
following properties are equivalent:
\begin{enumerate}[(i)]
\item $P$ embeds in
$[E]^{<\omega}$ as a join-subsemilattice for some set $E$;
\item   for every $x\in P$, $\varphi_{\triangle}(x)$ is finite.
\end{enumerate}
\end{proposition}
\begin{proof}
$(i)\Rightarrow (ii)$ Let $\varphi$ be an embedding from $P$ in
$[E]^{<\omega}$ which preserves joins.  Set $\mathcal F:= \varphi
(P)$. Apply   part $(c)$  of Lemma \ref{cl2} .
 $(ii)\Rightarrow (i)$ Set $E:= \triangle (J(P))$. We have  $\varphi_{\triangle}(x)\in [E]^{<\omega}$.  According to Fact \ref{representation}
 and Lemma \ref{lem01},  the map $\varphi_{\triangle}: P \rightarrow [E]^{<\omega}$ is an embedding preserving joins.\end{proof}

 \begin{corollary}\label{sierpinski} Let $\beta$ be a countable  order type. If a proper initial segment
 contains infinitely many non-principal  initial segments then no sierpinskisation  $P$ of  $\beta$ with $\omega$
 can embed  in $[\omega]^{<\omega}$ as a join-semilattice (whereas it embeds as a poset).
 \end{corollary}
 \begin{proof}
According to Proposition \ref{lem03}  it suffices to prove  that $P$ contains some $x$ for which $\varphi_{\Delta}(x)$ is infinite.

 Let $P$ be a sierpinskisation  of $\beta$  and  $\omega$. It is  obtained as the intersection of two linear orders $L$,   $L'$ on the same set  and having respectively order type  $\beta$  and  $\omega$.   We may suppose that the ground set is  $\N$ and $L'$ the natural order.

{\bf  Claim 1}
 A non-empty subset $I$ is a non-principal ideal  of
$P$ if and only if this is a non-principal initial segment of $L$.

{\bf Proof of  Claim 1} Suppose that  $I$ is  a non-principal
initial segment of $L$. Then, clearly, $I$ is an initial segment
of $P$. Let us check that $I$  is up-directed. Let $x, y\in I$;
since $I$ is non-principal in $L$, the set $A:= I\cap\uparrow_{L}
x \cap \uparrow_{L} y$ of upper-bounds of $x$ and $y$ w.r.t.  $L$
which belong to  $I$ is infinite; since $B:= \downarrow_{L'} x
\cup \downarrow_{L'} y$ is finite, $A\setminus B$ is non-empty. An
arbitrary element $z\in A\setminus B$ is an upper bound of $x,y$
in $I$ w.r.t. the poset  $P$ proving that $I$ is up-directed.
Since $I$ is infinite, $I$ cannot have a largest element in $P$,
hence $I$ is a non-principal ideal of $P$. Conversely, suppose
that $I$ is  a non-principal ideal of $P$.   Let us check that $I$
is an initial segment of $L$. Let $x\leq_{L} y$ with $y\in I$.
Since $I$ non-principal in $P$, $A:=\uparrow_{P} y\cap I$ is
infinite; since $B:= \downarrow_{L'} x \cup\downarrow_{L'} y$ is
finite, $A\setminus B$ is non-empty. An arbitrary element of
$A\setminus B$ is an upper bound of $x$ and $y$ in $I$ w.r.t. $P$.
It  follows that $x\in I$. If $I$ has a largest element w.r.t. $L$
then such an element must be maximal in $I$ w.r.t. $P$, and since
$I$ is an ideal, $I$ is a principal ideal, a contradiction.

{\bf Claim 2} Let  $x\in \N$. If there is a non-principal ideal of
$L$ which does not contain $x$, there is a maximal one, say $I_x$.
If $P$ is a join-semilattice, $I_x\in \Delta(P)$.

{\bf Proof of Claim 2} The first part follows from Zorn's Lemma. The second part follows from Claim 1 and
 Lemma  \ref{lem01}.

 {\bf Claim 3} If an initial segment $I$ of $\beta$ contains infinitely many non-principal initial segments
 then there is an infinite sequence $(x_n)_{n<\omega}$ of elements of $I$ such that the  $I_{x_n}$'s  are all distinct.

 {\bf Proof of Claim 3} With Ramsey's theorem obtain a sequence $(I_n)_{n<\omega}$ of non-principal initial segments which is either strictly
 increasing \index{increasing sequence} or strictly decreasing \index{decreasing sequence}.
 Separate two successive members by some element $x_n$ and apply the first part of Claim 2.

If we pick $x\in \N \setminus I$ then it follows from Claim 3 and the second part of Claim 2 that
$\varphi_{\Delta}(x)$ is infinite.
\end{proof}

\begin{example}\label{ex:ordinal} If $\alpha$ is a countably infinite order type distinct from  $\omega$, $\Omega(\alpha)$ is not embeddable  in $[\omega]^{<\omega}$ as a join-semilattice.
\end{example}
Indeed,  $ \Omega(\alpha)$ is a sierpinskisation of $ \omega\alpha$ and $\omega$.  And  if $\alpha$ is distinct from  $\omega$,
 $\alpha$ contains  some element which majorizes  infinitely many others.  Thus $\beta:= \omega\alpha$ satisfies
 the hypothesis of Corollary \ref{sierpinski}.

 Note that on an other hand, for every ordinal $\alpha\leq \omega$, there are  representatives of  $ \Omega(\alpha)$ which are  embeddable in $[\omega]^{<\omega}$ as  join-semilattices.

\begin{theorem}\label {finitesubsets}
Let $Q$ be a well-founded poset and let $\mathcal{F}\subseteq I_{<\omega}(Q)$. The following properties are equivalent:
\begin{itemize}
\item[$1)$] $\mathcal{F}$ has no infinite antichain;
\item[$2)$] $\mathcal{F}^{<\omega}$ is wqo;
\item[$3)$] $J(\mathcal{F}^{<\omega})$ is topologically scattered;
\item[$4)$] $\mathcal{F}^\cup$ is order-scattered;
\item[$5)$] $\mathfrak{P}(\omega)$ does not embed in $\mathcal{F}^\cup$;
\item[$6)$] $\lbrack \omega\rbrack^{<\omega}$  does not embed in $\mathcal{F}^{<\omega}$;
\item[$7)$] $\mathcal{F}^\cup$ is well-founded.
\end{itemize}
\end{theorem}

\begin{proof} We prove the following chain of implications:
$$1)\Longrightarrow 2)\Longrightarrow 3)\Longrightarrow
4)\Longrightarrow 5)\Longrightarrow 6)\Longrightarrow 7)
\Longrightarrow 1)$$ $1)\Longrightarrow 2)$. Since $Q$ is
well-founded then, as mentioned in $a)$ of Theorem
\ref{wellfounded}, $I_{<\omega}(Q)$ is well-founded.  It follows
first that $\mathcal{F}^{<\omega}$ is well-founded, hence from
Property $c)$ of Theorem  \ref{wellfounded},  every member of
$\mathcal{F}^{<\omega}$ is a finite join of join-irreducibles. Next,
as a subset of $\mathcal{F}^{<\omega}$, $\mathcal {F}$ is
well-founded, hence wqo according to  our hypothesis. The set of
join-irreducible members of $\mathcal{F}^{<\omega}$ is wqo as a
subset of $\mathcal {F}$.
 From Property $d)$  of Theorem  \ref{wellfounded}, $\mathcal{F}^{<\omega}$ is wqo \\
$2)\Longrightarrow 3)$. If $\mathcal{F}^{<\omega}$ is wqo then
$I(\mathcal{F}^{<\omega})$ is well-founded (cf. Property  ($b)$ of
Theorem \ref{wellfounded}). If follows that
$I(\mathcal{F}^{<\omega})$ is topologically scattered
(cf.\cite{misl});  hence all its subsets are topologically
scattered, in particular $J(\mathcal{F}^{<\omega})$. \\
$3)\Longrightarrow4)$. Suppose that $\mathcal{F}^\cup$ is not
ordered scatered. Let $f: \eta \rightarrow \mathcal{F}^\cup$ be an
embedding. For $r\in \eta$ set $\check f(r)=\bigcup \{f(r'):
r'<r\}$. Let $X:=\{\check f(r): r<\eta\}$. Clearly $X\subseteq
\mathcal{F}^\cup$. Furthermore $X$ contains no isolated point
(Indeed, since $\check f(r)=\bigcup \{\check f(r'): r'<r\}$, $\check
f(r)$ belongs to the topological closure of $\{\check f(r'):
r'<r\}$). Hence $\mathcal{F}^\cup$ is not topologically scatered.\\
$4)\Longrightarrow 5)$. Suppose that $\mathfrak{P}(\omega)$ embeds
in $\mathcal{F}^\cup$. Since $\eta\leq \mathfrak{P}(\omega)$, we
have
$\eta\leq \mathcal{F}^\cup$.\\
$5)\Longrightarrow 6)$. Suppose that $[\omega]^{<\omega}$ embeds in
$\mathcal{F}^{<\omega}$, then $J([\omega]^{<\omega})$ embeds in
$J(\mathcal{F}^{<\omega})$. Lemma \ref{cl2} assures that
$J(\mathcal{F}^{<\omega})$ is isomorphic to $\mathcal{F}^{\cup}$. In
the other hand $J([\omega]^{<\omega})$ is isomorphic to
$\mathfrak{P}(\omega)$. Hence $\mathfrak{P}(\omega)$ embeds in
$\mathcal{F}^{\cup}$.\\
$6)\Longrightarrow 7)$. Suppose $\mathcal{F}^{\cup}$ not
well-founded. Since $Q$ is well-founded, $a)$ of Theorem
\ref{wellfounded} assures $I_{<\omega}(Q)$ well-founded, but
$\mathcal{F}^{<\omega}\subseteq I_{<\omega}(Q)$, hence
$\mathcal{F}^{<\omega}$ is well-founded. Furthermore, since
$I_{<\omega}(Q)$ is closed under finite unions, we have
$\mathcal{F}^{<\omega}\subseteq I_{<\omega}(Q)$, Proposition \ref
{w-f} implies that $\underline\Omega (\omega^{*})$ does not embed in
$\mathcal{F}^{<\omega}$.
From Theorem \ref{thm4}, we have $\mathcal{F}^{<\omega}$ not well-founded.\\
$7)\Longrightarrow 1)$. Clearly, $\mathcal F$ is well-founded. If
$F_{0}, \dots, F_{n}\dots $ is an infinite antichain of members of
$\mathcal{F}$, define $f(i,j):[\omega]^{2}\rightarrow Q$, choosing
$f(i,j)$ arbitrary in $Max(F_{i})\setminus F_{j}$. Divide
$[\omega]^3$ into $R_{1}:=\{(i,j,k)\in [\omega]^3: f(i,j)=f(i,k)\}$
and $R_{2}:=[\omega]^3\setminus R_{1}$.
 From Ramsey's theorem,  cf. \cite {rams}, there is some infinite subset $X$ of $\omega$ such that $[X]^{3}$ is included in $R_{1}$ or in $R_{2}$.
The inclusion in $R_{2}$ is impossible since $\{f(i,j ): j<\omega
\}$, being included in $Max(F_{i})$,  is finite for every $i$. For
each $i\in X$, set $G_{i}:= \bigcup\{F_{j}: i\leq j\in X\}$. This
defines an $\omega^*$-chain in $\mathcal{F}^\cup$.
\end{proof}

\begin{remark} If $\mathcal{F}^{<\omega}$ is closed under finite
intersections then equivalence between $(3)$ and $(4)$ follows from
Mislove's Theorem mentioned in \cite{misl}.
\end{remark}
Theorem \ref {finitesubsets} above was obtained by the second author and M.Sobrani in the special case where $Q$ is an antichain \cite {pouzet, sobrani} .

 \begin{corollary} \label{provisoire}
If  $P$ is a join-subsemilattice of a join-semilattice of the form $[\omega]^{<\omega}$,
or more generally of the form $I_{<\omega} (Q)$ where $Q$ is some well-founded poset,
then $J(P)$  is well-founded if and only if $P$ has no infinite antichain.
\end{corollary}

{\bf Remark.} If, in Theorem \ref {finitesubsets} above, we suppose that $\mathcal  F$ is well-founded instead of $Q$,
all implications in the above chain hold, except $6)\Rightarrow 7)$.  A counterexample is provided by $Q:= \omega \oplus \omega^*$,
the direct sum of the chains $\omega$ and
$\omega^*$, and
$\mathcal F$, the image of $\underline \Omega (\omega^*)$ via a natural embedding.

\subsection{ Proof of  Theorem \ref {thmwf}}
 $(i)\Rightarrow (ii)$ Suppose that $(i)$ holds. Set $Q:= J(P)$. Since $P$ contains no infinite antichain,
 $P$ embeds as a join-subsemilattice in $ I_{<\omega}(Q)$ (Corollary \ref {corthm16}). From $b)$ of Theorem \ref{wellfounded} $Q$ is well-founded.
 Since $P$ has no infinite antichain, it has no infinite independent set.

  $(ii)\Rightarrow (i)$ Suppose that $(ii)$ holds. Since $Q$ is well-founded, then from $a)$ of Theorem \ref{wellfounded}, $I_{<\omega}(Q)$ is well-founded.
  Since $P$ embeds in $I_{<\omega}(Q)$, $P$ is well-founded. From our hypothesis, $P$ contains no infinite independent set.  According to  implication $(iii)\Rightarrow (i)$ of Theorem \ref{tm2.1} , it does not embed $[\omega]^{<\omega}$. From implication $6)\Rightarrow 1)$ of Theorem \ref {finitesubsets}, it has no infinite antichain. \endproof

\end{document}